\documentclass{article}
\usepackage[english]{babel}
\usepackage[latin1]{inputenc}
\usepackage{latexsym}
\usepackage{amssymb,amsmath,dsfont}
\usepackage{subfigure}
\usepackage{graphicx,float}
\usepackage{color}
\usepackage{url}
\usepackage{authblk}

\usepackage{array}
\usepackage[table]{xcolor}
\usepackage{amsthm}
\usepackage{a4}
\usepackage{fullpage}
\usepackage[unicode]{hyperref}
\hypersetup{
    bookmarks=true,         % show bookmarks bar?
    unicode=true,          % non-Latin characters in Acrobat?s bookmarks
    pdftoolbar=true,        % show Acrobat?s toolbar?
    pdfmenubar=true,        % show Acrobat?s menu?
    pdffitwindow=true,      % page fit to window when opened
    pdftitle={FKA},    % title
    pdfauthor={C. Bourdarias, M. Ersoy and S. Gerbi},     % author
    pdfsubject={},   % subject of the document
    pdfcreator={},   % creator of the document
    pdfproducer={}, % producer of the document
    pdfkeywords={}, % list of keywords
    pdfnewwindow=true,      % links in new window
    colorlinks=true,       % false: boxed links; true: colored links
    linkcolor=blue,  
    citecolor=blue,        % color of links to bibliography
    filecolor=magenta,      % color of file links
    urlcolor=cyan           % color of external links
}
\newtheorem{theorem}{Theorem}[section]

\newtheorem{lemma}{Lemma}[section]

\newtheorem{definition}{Definition}[section]

\newtheorem{remark}{Remark}[section]
\newtheorem*{ack}{Acknowledgments}

\newcommand{\R}{\mathbb{R}}
\makeatletter
\newcommand{\bigint}{\@ifnextchar_\@bigintsub\@bigintnosub}
\def\@bigintsub_#1{\def\@int@subscript{#1}\@ifnextchar^\@bigintsubsup\@bigintsubnosup}
\def\@bigintsubsup^#1{\mathop{\text{\large$\int_{\text{\normalsize$\scriptstyle\@int@subscript$}}^{\text{\normalsize$\scriptstyle#1$}}$}}\nolimits}
\def\@bigintsubnosup{\mathop{\text{\large$\int_{\text{\normalsize$\scriptstyle\@int@subscript$}}$}}\nolimits}
\def\@bigintnosub{\@ifnextchar^\@bigintnosubsup\@bigintnosubnosup}
\def\@bigintnosubsup^#1{\mathop{\text{\large$\int^{\text{\normalsize$\scriptstyle#1$}}$}}\nolimits}
\def\@bigintnosubnosup{\mathop{\text{\large$\int$}}\nolimits}
\makeatother
 \title{Asymptotic stability and blow up for a semilinear damped wave equation with dynamic boundary conditions.}
\author[1]{{\sc{S. Gerbi}}\thanks{stephane.gerbi@univ-savoie.fr}}
\author[2]{{\sc{B. Said-Houari}}\thanks{saidhouarib@yahoo.fr}}
\affil[1]{\small{Laboratoire de Math\'ematiques, UMR 5127 - CNRS and Universit\'e de Savoie, \par73376 Le Bourget-du-Lac Cedex, France.}}
\affil[2]{\small{Division of Mathematical and Computer Sciences and Engineering, \par King Abdullah University of Science and Technology (KAUST), Thuwal, Saudi Arabia}}
\begin{document}
\date{}
\maketitle
\begin{abstract}
In this paper we consider a multi-dimensional wave equation with dynamic boundary conditions, related to the Kelvin-Voigt damping. Global existence and asymptotic stability of solutions starting in a stable set are proved. Blow up for solutions
of the problem with linear dynamic boundary conditions with initial data in the unstable set is also obtained.
\end{abstract}
\section{Introduction}
In this paper we consider the following semilinear damped wave equation with dynamic boundary conditions:
%%%%%%%%%%%%%%%%%%%%%%%%%%%%%%%%%%%%%%%%%%%
% L'EDP hyperbolique degeneree
%%%%%%%%%%%%%%%%%%%%%%%%%%%%%%%%%%%%%%%%%%%%
\begin{equation}\label{ondes}
\left\{
\begin{array}{ll}
%%% equation
\hspace*{-0.0cm}u_{tt}-\Delta u-\alpha \Delta u_{t}=\vert u\vert^{p-2}u, & x\in \Omega ,\ t>0 \\[0.1cm]
%%% condtions aux limites de Dirichlet
\hspace*{-0.0cm}u( x,t) =0, &  x\in \Gamma_{0},\ t>0  \\[0.1cm]
%%% condtions aux limites mixtes nonlineaire
\hspace*{-0.0cm}u_{tt}(x,t) =- a \left[\displaystyle \frac{\partial u}{\partial \nu }(x,t) +
\frac{\alpha \partial u_{t}}{\partial \nu }(x,t) +r\vert u_{t}\vert^{m-2}u_{t}( x,t) \right]
&  x\in \Gamma_{1},\ t>0 \\[0.1cm]
%%% conditions initiales
\hspace*{-0.25cm}u( x,0) =u_{0}(x), \; u_{t}( x,0) =u_{1}(x) & x\in \Omega ,
\end{array}
\right.
\end{equation}
%%%%%%%%%%%%%%%%%%%%%%%%%%%%%%%%%%%%%%%%%%%%%%%%%%%%%%%%%%%%%%%%%%%%%%%%%%%%%%%%%%%%%%%%%%%%%%%%%%%
where $u=u(x,t)\,,\, t \geq 0\,,\, x\in \Omega\,,\,\Delta$ denotes the Laplacian operator with res\-pect to the $x$ variable,  $\Omega$ is a regular and bounded domain  of $\mathbb{R}^{N}\,,\,( N\geq 1)$, $\partial\Omega~=~\Gamma_{0}~\cup~\Gamma_{1}$, $mes(\Gamma_{0}) >0,$
$\Gamma_{0}\cap \Gamma_{1}=\varnothing $ and $\displaystyle \frac{\partial }{\partial \nu }$
denotes the unit outer normal derivative, $m \geq 2\,,\,a\,,\,\alpha \mbox{ and } r$ are
positive constants, $p>2 $ and $u_{0}\,,\,u_{1}$ are given functions. For the sake of simplicity, in this paper we consider the problem (\ref{ondes}) where we have set
$a = 1$. 

From the mathematical point of view, these problems do not neglect  acceleration terms on the boundary.
Such type of boundary conditions are usually called \textit{dynamic boundary conditions}.
They are not only important from the theoretical point of view but also arise in numerous practical problems.
For instance in one space dimension, the problem (\ref{ondes}) can modelize the dynamic evolution of a viscoelastic
rod that is fixed at one end and has a tip mass attached to its free end.
The dynamic boundary conditions represent  the Newton's law for the attached mass, (see \cite{BST64,AKS96, CM98} for more details).
In the two dimension space, as showed in \cite{G06} and in the references therein, these  boundary conditions arise  when we consider the transverse motion of a flexible membrane
$\Omega$ whose boundary may be affected by the vibrations only in a region.
Also some dynamic boundary conditions as in problem (\ref{ondes}) appear when we assume that $\Omega$ is an exterior domain of $\mathbb{R}^3$ in  which homogeneous fluid is at rest except for sound waves. Each point of the boundary is subjected to small normal
displacements into the obstacle  (see \cite{B76} for more details).
This type of dynamic boundary conditions are known as acoustic boundary  conditions.
More  results on the wave equations with acoustic boundary conditions can be found in \cite{FG07_1}.

Before state and prove our results, let us first recall some works related to the problem we address.

Among the early results dealing with this type of boundary conditions are those of Grobbelaar-Van Dalsen \cite{GV94,GV96}
in which the author has made contributions to this field.

In \cite{GV94} the author introduced a model which describes the damped longitudinal vibrations of a homogeneous flexible
horizontal rod of length $L$ when the end $x = 0$ is rigidly fixed while the other end $x=L$ is free to move with an attached load.
This yields to a following systems of partial differential equations:
\begin{equation}\label{Vad_1}
\left\{
\begin{array}{ll}
u_{tt}-u_{xx}-u_{txx}=0, & x\in (0,L),\,t>0 \\[0.1cm]
u(0,t)=u_{t}(0,t)=0, & t>0 \\[0.1cm]
u_{tt}(L,t)=-\left[ u_{x}+u_{tx}\right] (L,t), & t>0 \\[0.1cm]
u\left( x,0\right) =u_{0}\left( x\right) ,u_{t}\left( x,0\right)
=v_{0}\left( x\right)  & x\in (0,L)\\
u\left( L,0\right) =\eta ,\qquad u_{t}\left( L,0\right) =\mu \ .
\end{array}
\right.
\end{equation}
By rewriting problem (\ref{Vad_1}) within the framework of the abstract theories of the so-called $B$-evolution theory, an
existence of a unique solution in the strong sense has been shown. An exponential decay result was also shown in \cite{GV96}
for a problem related to (\ref{Vad_1}), which describe the
weakly damped vibrations of an extensible beam. See \cite{GV96} for more details.

Subsequently, Zang and Hu \cite{ZH07}, considered the problem
\begin{equation}\label{Zang_Hu}
\left\{
\begin{array}{ll}
u_{tt}-p\left( u_{x}\right) _{xt}-q\left( u_{x}\right) _{x}=0,&x\in \left(
0,1\right) ,\,t>0, \vspace{0.2cm}\\
u\left( 0,t\right) =0,&t\geq 0\vspace{0.2cm}\\
\left( p\left( u_{x}\right) _{t}+q\left( u_{x}\right)
\left( 1,t\right) +ku_{tt}\left( 1,t\right) \right) =0,&t\geq 0, \vspace{0.2cm}\\
u\left( x,0\right) =u_{0}\left( x\right) ,\qquad u_{t}\left( x,0\right)
=u_{1}\left( x\right) ,&x\in \left( 0,1\right) \ .
\end{array}%
\right.
\end{equation}
By using the Nakao inequality, and under appropriate conditions on $p$ and $q$, they established both an exponential
and polynomial decay rates for the energy depending on the form of the terms $p$ and $q$.

It is clear that in the absence of the source term $\vert u\vert^{p-2}u$ and for $r=0$, problem (\ref{Vad_1})
 is the one dimensional model of (\ref{ondes}). Similarly, in the case where the source term vanishes
  identically and for $r\neq 0 \,,\,  m = 2 \mbox{ and } p = 2$, Pellicer and Sol{\`a}-Morales \cite{PS04}
considered the one dimensional problem as an alternative model for the classical spring-mass damper system,
and by using the dominant eigenvalues method, they showed that the large time behavior of the solutions of problem (\ref{ondes})
 is the same as for a classical spring-mass damper ODE, namely:
\begin{equation}\label{ODE}
m_1 u''(t)+ d_1 u'(t)+k_1 u(t)=0,
\end{equation}
when $a$ tends to zero, where the parameters $m_1 \,,\, d_1 \mbox{ and } k_1$ are determined from the values of the spring-mass damper system.

Thus, the asymptotic stability of the model for small values of $a$ has been determined as a consequence of this limit. But they did not obtain any rate of convergence.
This result was followed by recent works \cite{P08,PS08}. In \cite{PS08}, a
continuous model for a spring-mass-damper system has been treated, where possible differences
 in the internal deformation of the spring are considered. More precisely, they investigated the following problem
 \begin{equation}\label{Pellicer_Mora_2}
    \left\{
\begin{array}{ll}
%%% equation
u_{tt}- u_{xx}-\alpha u_{txx} =0,&x\in (0,1),\,t>0  \\[0.1cm]
%%% condtions aux limites de Dirichlet
u( 0,t) =0,&t>0  \\[0.1cm]
%%% condtions aux limites mixtes nonlineaire
u_{tt}(1,t) =- \varepsilon\left[u_{x} +\alpha u_{tx}+ru_{t}\right](1,t)
,&t>0  \ .\\[0.1cm]
%%% conditions initiales
\end{array}
\right.
 \end{equation}
 By using the spectral analysis approach, they showed that for different values
of the parameter $\alpha$, the limit behaviors are very different from the classical
ODE (\ref{ODE}).
While in \cite{P08} the author considered a one dimensional nonlocal nonlinear strongly damped wave equation with dynamical boundary conditions.
In other words, they looked to the following problem:
\begin{equation}\label{spring-mass}
\left\{
\begin{array}{ll}
%%% equation
u_{tt}- u_{xx}-\alpha u_{txx}+\varepsilon f\left( u(1,t),\frac{u_{t}(1,t)}{\sqrt{\varepsilon}}\right) =0,\  \\[0.1cm]
%%% condtions aux limites de Dirichlet
u( 0,t) =0,  \\[0.1cm]
%%% condtions aux limites mixtes nonlineaire
u_{tt}(1,t) =- \varepsilon\left[u_{x} +\alpha u_{tx}+ru_{t}\right](1,t)-\varepsilon f\left( u(1,t),\frac{u_{t}(1,t)}{\sqrt{\varepsilon}}\right)
,\ \\[0.1cm]
%%% conditions initiales
\end{array}
\right.
\end{equation}
%%%%%%%%%%%%%%%%%%%%%%%%%%%%%%%%%%%%%%%%%%%%%%%%%%%%%%%%%%%%%%%%%%%%%%%%%%%%%%%%%%%%%%%%%%%%
with $x\in(0,1),\,t>0,\,r,\alpha>0$ and $\varepsilon\geq0$.
The above system modelises a spring-mass-damper system, where the term $\varepsilon f\left( u(1,t),\frac{u_{t}(1,t)}
{\sqrt{\varepsilon}} \right)$ represents a control acceleration at $x=1$.
By using the invariant manifold theory, the authors proved that for small values of the parameter $\varepsilon$,
the solutions of  (\ref{spring-mass}) are attracted to a two dimensional invariant manifold.
See \cite{PS08}, for further details.

We recall that the presence of the strong damping term $-\Delta u_{t}$ in the problem (\ref{ondes}) makes the problem different
from that considered in \cite{GT94} and widely studied in the literature \cite{V99,T98,T99,GS06,TV05} for instance.
For this reason less results were known for the wave equation with a strong damping and many problems remained unsolved.
Especially the blow-up of solutions in the presence of a strong damping and a nonlinear boundary damping at the same time
is still an open problem. In \cite{GS08}, the present authors showed that the solution of (\ref{ondes})
is unbounded and grows up exponentially when time goes to infinity if the initial data are large enough.

Recently, Gazzola and Squassina \cite{GS06} studied the global solution and the finite time blow-up for a damped semilinear wave equation with Dirichlet boundary conditions by a careful study
of the stationary solutions and their stability using the Nehari manifold and a mountain pass energy level of the initial condition.

The main difficulty of the problem considered is related to the non ordinary boundary conditions defined on $\Gamma_1$.
Very little attention has been paid to this type of boundary conditions. We mention only a few particular results in the one dimensional space and for a linear damping i.e. $(m=2)$ \cite{GV99,PS04,DL02,K92}.

A problem related to (\ref{ondes}) is the following:
\begin{eqnarray}\label{Dron}
u_{tt}-\Delta u + g(u_{t}) &=&f \hspace*{1.5cm} \text{ in } \Omega \times ( 0,T)  \nonumber \\
\frac{\partial u}{\partial \nu }+ K(u) u_{tt}+ h(u_{t}) &=&0,\hspace*{1.5cm}
\text{ on }\partial \Omega \times (0,T)\\
u(x,0) &=& u_{0}(x) \hspace*{1cm}\text{ in }\Omega \nonumber \\
u_{t}(x,0)&=&u_{1}(x) \hspace*{1cm}\text{ in }\Omega  \nonumber
\end{eqnarray}
where $f=f(x,t)$ and the boundary term $h( u_{t}) =\vert u_{t}\vert^{\rho }u_{t}$ arises when one studies flows of gas in a channel with porous walls. The term $u_{tt}$ on the boundary appears from the internal
forces, and the nonlinearity $K(u) u_{tt}$ on the boundary represents the internal
 forces when the density of the medium depends on the displacement.
This problem has been studied in \cite{DL02}, in the one dimensional case and in
 \cite{DLS98} for $N$-dimensional with $N\geq1$. By using the Faedo-Galerkin
 approximations and a compactness argument, they proved the global existence
 of the solution. Also, the exponential decay of the total energy of problem (\ref{Dron}) has been proved
under the condition $f=0$.

Most of the above mentioned papers only treat particular cases of problem (\ref{ondes}). The aim of our previous paper \cite{GS08} and of
this paper is to apply known methods in order to investigate the more general problem (\ref{ondes}).

Recently, the present authors studied problem (\ref{ondes}) in \cite{GS08}. A local existence result was obtained by combining the Faedo-Galerkin method with the contraction mapping theorem. Concerning the asymptotic behavior, the authors showed that the solution of such problem is unbounded and grows up
exponentially when time goes to infinity if the initial data are large enough and the damping term is nonlinear (i.e. $m>2$).

As we have said before, our problem (\ref{ondes}) can be seen as a model which describe the interaction between an elastic medium and a rigid mass. So, it
seems more convenient to recall some results related to the interaction of an elastic medium with rigid mass. By using the classical semigroup theory, Littman and Markus \cite{LM88} established a uniqueness result for a particular Euler-Bernoulli beam rigid
body structure. They also proved the asymptotic stability of the structure by using the feedback boundary damping. In \cite{LL98} the authors considered the Euler-Bernoulli beam equation which describes the dynamics of clamped
elastic beam in which one segment of the beam is made with viscoelastic material and the other of elastic material. By combining the frequency domain method with the multiplier technique, they proved the exponential decay
for the transversal motion but not for the longitudinal motion of the model, when the Kelvin-Voigt damping is distributed only on a subinterval of the domain.
In relation with this point, see also the work by Chen et \textit{al.} \cite{CLL98} concerning the Euler-Bernoulli beam equation with the global or local Kelvin-Voigt damping. Also models of vibrating strings with local viscoelasticity and
Boltzmann damping, instead of the Kelvin-Voigt one, were considered in \cite{LL02} and an exponential energy decay rate was established. Recently, Grobbelaar-Van Dalsen \cite{G03} considered an extensible thermo-elastic beam
which is hanged at one end with rigid body attached to its free end, i.e. one dimensional hybrid thermoelastic structure, and showed that the method used in \cite{O97} is still valid to establish an uniform stabilization of the system.
Concerning the controllability of the hybrid system we refer to the work by Castro and Zuazua \cite{CZ98}, in which they considered flexible beams connected by point mass and the model takes account of the rotational inertia.

In this paper we consider the problem (\ref{ondes}) and we will show that if the initial data are in the ``stable set'', the solution continues to live there forever.
In addition, we will prove that the presence of the strong damping forces the solution to go to zero uniformly and with an exponential decay rate, even if the boundary damping is nonlinear i.e. $m > 2$.
To obtain our results we combine the potential well method with the energy method. We will also proved
that in the absence of the nonlinearity in the boundary term (that is, in the case where $m=2$), the solution blows up in finite time.

Let us now give a short summary of the content of this paper.
In section 2, after having stated the local existence and uniqueness theorem obtained by the authors in \cite{GS08},
we will prove that if the initial data are in the stable manifold, the solution continues to live there forever
and so we will prove the global existence and the exponential decay of the solution.

In section 3, we prove the blow up result of the problem (\ref{ondes}), in the case of a linear boundary damping (that is, $m=2$),
in spite of the presence of the strong damping term $\triangle u_{t}$.
The technique we use follows closely the method used in \cite{GS06}, which is based on the concavity argument due to Levine
\cite{L74_1}.
%%%%%%%%%%%%% NEW %%%%%%%%%%%
Let us mention, that despite the methods used here are well-known tools to
prove the global existence, exponential decay and blow of solution, therefore, the main novelty of the work presented in this paper is  that we will use these techniques
to study the asymptotic behavior of the semilinear damped wave equation with {\it dynamic boundary conditions}.
To our knowledge, this has not been done before and this is the first paper dealing with the asymptotic behavior of such problem.
%%%%%%%%%%%%% FIN NEW %%%%%%%%%%%
%%%%%%%%%%%%%%%%%%%%%%%%%%%%%%%%%%%%%%%%%%%%%
\section{Asymptotic stability}
%%%%%%%%%%%%%%%%%%%%%%%%%%%%%%%%%%%%%%%%%%%%%
In this section,  we will first recall the local existence and the uniqueness result of
the solution of the problem (\ref{ondes}) proved in \cite{GS08}.
Then we state and prove the global existence and exponential decay of the solution of problem (\ref{ondes}).
In order to do this, a suitable choice of the Lyapunov functional will be made.\\

Let us first present some material that we shall use later in this paper.
We denote
$$
H_{\Gamma_{0}}^{1}(\Omega) =\left\{u \in H^1(\Omega) /\ u_{\Gamma_{0}} = 0\right\} .
$$
By $( .,.) $ we denote the scalar product in $L^{2}( \Omega)$  i.e. $(u,v)(t) = \displaystyle \int_{\Omega} u(x,t) v(x,t) dx$. Also we mean by $\Vert .\Vert_{q}$ the $L^{q}(\Omega) $ norm for $1 \leq q \leq \infty$, and by
$\Vert .\Vert_{q,\Gamma_{1}}$ the $L^{q}(\Gamma_{1}) $ norm.

Let us denote for $v \in H_{\Gamma_{0}}^{1}(\Omega)$
\begin{equation}\label{normstar}
\Vert v\Vert_{\ast }^{2}=\alpha \Vert v\Vert _{2,\Gamma_{1}}^{2}+r\Vert \nabla v\Vert _{2}^{2}
\end{equation}
Let $T>0$ be a real number and $X$ a Banach space endowed with norm $\Vert .\Vert _{X}$.
$L^{p}(0,T;X) ,\ 1 \leq p < \infty$ denotes the space of functions $f$ which are $L^{p}$
over $\left( 0,T\right) $ with values in $X$, which are measurable and  $\Vert f \Vert_{X} \in L^{p} \left(0,T\right)$.
This space is a Banach space endowed with the norm
$$
\Vert f\Vert_{L^{p}\left( 0,T;X\right) }=
\left(\int_{0}^{T}\Vert f\Vert_{X}^{p} dt\right)^{1/p}.
$$
$L^{\infty}\left( 0,T;X\right) $ denotes the space of functions $f:\left]0,T\right[ \rightarrow X$ which are measurable and
$\Vert f\Vert_{X}\in L^{\infty }\left( 0,T\right) $.
This space is a Banach space endowed with the norm:
$$
\Vert f\Vert_{L^{\infty}(0,T;X)}=\mbox{ess}\sup_{0<t<T}\Vert f\Vert_{X}.
$$
We recall that if $X$ and $Y$ are two Banach spaces such that $X\hookrightarrow Y$ (continuous embedding), then
$$
L^{p}\left( 0,T;X\right) \hookrightarrow L^{p}\left( 0,T;Y\right) , \ 1 \leq p\leq \infty .
$$
 We define the critical Sobolev exponent for the trace functional space by:
\begin{equation}\label{qbar}
\bar{q}=\left\{
\begin{array}{cl}
\displaystyle \frac{2 (N-1) }{N-2} \,,&\mbox{ if } N \geq 3 \\
+\infty \,,&\mbox{ if } N=1, 2.
\end{array}
\right.
\end{equation}
Let us define the space $Y_T $ as:
\begin{eqnarray*}
Y_T =  &\biggl\{&(v,v_t): v\in C \Bigl( [ 0,T],H_{\Gamma_{0}}^{1}(\Omega) \Bigl) \cap C^{1}\Bigl( [ 0,T] ,L^{2}(\Omega) \Bigl), \\
&&v_{t} \in L^{2}\Bigl(0,T;H_{\Gamma_{0}}^{1}(\Omega)\Bigl) \cap L^{m}\Bigl( (0,T) \times \Gamma_{1}\Bigl)\biggr\}
\end{eqnarray*}
endowed with the norm:
$$
\Vert (v,v_{t}) \Vert_{Y_T}^2=\max_{0 \leq t \leq T}\Bigl[ \Vert v_{t}\Vert_{2}^2 +
\Vert \nabla v \Vert_{2}^2\Bigl] + \Vert v_{t} \Vert^2_{L^{m}\bigl( \left( 0,T\right) \times \Gamma_{1}\bigl) } +
\int_0^T \Vert \nabla v_t(s) \Vert_2^2 \; ds \;  .
%%%%%%%%%%%%%%% NEW %%%%%%%%%%%%%%%%%
$$
For $m \leq \bar{q}$,  from Poincar\'e's inequality, the continuity of the trace operator on $\Gamma_{1}$ and 
Sobolev imbedding this norm is equivalent to:
\begin{equation}\label{norm}
\Vert u \Vert = \max_{0 \leq t \leq T}\left[\Vert \nabla u \Vert_2  + \Vert u_t \Vert_2 \right].
\end{equation}

%%%%%%%%%%%%%%%%%%%%%%%%%%%%%%%%
In this work, we will deal with the weak solution of the problem (\ref{ondes}),
consequently, we use the same definition as in \cite{GS08}.
\begin{definition} \label{generalised}
A function $u(x,t) $ defined on $\Omega \times [0,T]$, such that
\begin{eqnarray*}
u &\in &L^{\infty}\left( 0,T;H_{\Gamma_{0}}^{1}(\Omega) \right) \ , \\
u_{t} &\in &L^{2}\left( 0,T;H_{\Gamma_{0}}^{1}(\Omega) \right) \cap L^{m}\left( ( 0,T) \times \Gamma_{1}\right) \ , \\
u_{t} &\in &L^{\infty}\left( 0,T;H_{\Gamma_{0}}^{1}(\Omega)\right) \cap L^{\infty}\left( 0,T;L^{2}(\Gamma_{1})\right)\ , \\
u_{tt} &\in &L^{\infty}\left( 0,T;L^{2}(\Omega)\right) \cap L^{\infty}\left( 0,T;L^{2}(\Gamma_{1}) \right) \ , \\
u(x,0) & = &u_{0}(x)\,,\\
u_{t}(x,0) &= & u_{1}(x) \,,
\end{eqnarray*}%
is a generalized solution to the problem (\ref{ondes}) if for any function
$\omega \in H_{\Gamma_{0}}^{1}(\Omega) \cap L^{m}(\Gamma_{1}) $ and $\varphi \in C^{1}(0,T) $ with
$\varphi(T) =0$, we have the following identity:
\begin{equation*}
\begin{array}{lll}
\hspace*{-0.2cm}\displaystyle  \int_{0}^{T} (|u|^{p-2}u,w)(t) \, \varphi(t)\, dt &=& \displaystyle \int_{0}^{T}
\Bigl[ (u_{tt},w)(t) +(\nabla u,\nabla w)(t) + \alpha ( \nabla u_{t},\nabla w)(t)\Bigl] \, \varphi(t)\, dt \\
&+&\displaystyle \int_{0}^{T}\varphi(t) \int_{\Gamma_{1}}
\Bigl[u_{tt}(t) +r \vert u_{t}(t) \vert^{m-2}u_{t}(t) \Bigl] w \, d\sigma \, dt .
\end{array}%
\end{equation*}
\end{definition}
%%%%%%%%%%%%%%%%%%%%%%%%%%%%%%%%%
\begin{theorem}\label{existence}  {\rm{ \cite{GS08}}}
Let $2\leq p\leq \bar{q}$ and $\max\left( 2,\displaystyle \frac{\bar{q}}{\bar{q}+1-p} \right) \leq m \leq \bar{q}$. \\
Then given $u_{0}\in H_{\Gamma_{0}}^{1}(\Omega) $ and $u_{1}\in L^{2}(\Omega) $, there exists $T > 0$ and a unique
solution $u$ of the problem (\ref{ondes}) on $[0,T) $ such that
\begin{eqnarray*}
u &\in &C \Bigl( [ 0,T],H_{\Gamma_{0}}^{1}(\Omega) \Bigl) \cap C^{1}\Bigl( [ 0,T] ,L^{2}(\Omega) \Bigl), \\
u_{t} &\in &L^{2}\Bigl(0,T;H_{\Gamma_{0}}^{1}(\Omega)\Bigl) \cap L^{m}\Bigl( (0,T) \times \Gamma_{1}\Bigl)
\end{eqnarray*}
\end{theorem}
%%%%%%%%%%%%%%%%%%%%%%%%%%%%%%%%%
We proved this theorem by using the Faedo-Galerkin approximations and the well-known contraction mapping theorem.

%%%%%%%%%%%%%%%%%%%%%%%%%%%%%%%%%%%%%%%%%%%
% Defintion Tmax
%%%%%%%%%%%%%%%%%%%%%%%%%%%%%%%%%%%%%%%%%%%
\begin{definition} \label{Tmax}
Let $2\leq p\leq \bar{q}$, $\max\left( 2,\displaystyle \frac{\bar{q}}{\bar{q}+1-p} \right) \leq m \leq \bar{q}$,
$u_{0}\in H_{\Gamma_{0}}^{1}(\Omega) $ and $u_{1}\in L^{2}(\Omega) $. We denote by $u$ the solution of (\ref{ondes}).
We define:
$$
T_{max} = \sup\Bigl\{ T > 0 \,,\, u = u(t) \; exists \; on \; [0,T]\Bigr\}
$$
Since the solution $u \in Y_T$ (the solution is ``regular enough''), from the definition of the norm given by (\ref{norm}), let us recall that if
$T_{max} < \infty$, then
$$
 \lim_{\underset {t < T_{max}} {t \rightarrow T_{max}}} \Vert \nabla u \Vert_2  + \Vert u_t \Vert_2 = + \infty.
$$
If $T_{max} < \infty$, we say that the solution of (\ref{ondes}) blows up and that $T_{max}$ is the blow up time.\\
If $T_{max} = \infty$, we say that the solution of (\ref{ondes})  is global.
\end{definition}

In order to study the blow up phenomenon or the global existence of the solution of (\ref{ondes}), and following \cite{GS06},  
we define the functions $ I,J:H_{\Gamma_{0}}^{1}(\Omega)\mapsto \R$ by:
%%%%%%%%
\begin{eqnarray*}
I(u) &=&\Vert \nabla u \Vert_{2}^{2} - \Vert u \Vert_{p}^{p}  \ , \\
J(u) &=&\frac{1}{2} \Vert \nabla u \Vert_{2}^{2}-\frac{1}{p} \Vert u \Vert_{p}^{p} \ .
\end{eqnarray*}
For a given function $u \in H_{\Gamma_{0}}^{1}(\Omega)$, when we will use the evaluation of the above functions at a time $0 \leq t < T_{max}$,  for the sake of simplicity,  we will write:
%%%%%%%%
\begin{eqnarray}
I( t) & =&\Vert \nabla u(t) \Vert_{2}^{2} - \Vert u(t) \Vert_{p}^{p} \ ,\label{Energy_I} \\
J(t) &  =&\frac{1}{2} \Vert \nabla u(t) \Vert_{2}^{2}-\frac{1}{p} \Vert u(t) \Vert_{p}^{p} \ . \label{Energy_J} 
\end{eqnarray}
We then define the energy of a solution $u$ of (\ref{ondes}) as:
\begin{equation}
E(t)= J(t) +\frac{1}{2}\Vert u_{t}(t)\Vert_{2}^{2}+\frac{1}{2}\Vert u_{t}(t)\Vert_{2,\Gamma_{1}}^{2}\, \quad  \forall \, 0 \leq t< T_{max} \label{Energy_E}
\end{equation}
As in \cite{GS08}, multiplying the first equation in (\ref{ondes}) by $u_{t}$  and integrating over $\Omega$ and with respect to $t$,
 we obtain the following  energy identity :
\begin{equation} \label{derivE}
E(t) - E(s) =- \int_{s}^{t} \Vert u(\tau) \Vert_{\ast}^{2}d\tau,\quad \forall \, 0 \leq s \leq t < T_{max}.
\end{equation}
Thus the function $E$ is decreasing along the trajectories.

As in \cite{PS75}, the potential well depth is defined as:
\begin{equation}\label{potentialwell}
d=\inf_{u\in H_{\Gamma_{0}}^{1}(\Omega)\backslash \{0\}} \max_{\lambda \geq 0}J(\lambda u) .
\end{equation}
We can now define the so called ``Nehari manifold'' as follows:
\begin{equation*}
\mathcal{N}=\left\{ u\in H_{\Gamma_{0}}^{1}(\Omega) \backslash \{0\} ; \; I(u) =0\right\} .
\end{equation*}
$\mathcal{N}$ separates the two unbounded sets:
\begin{equation*}
\mathcal{N}^{+}= \left\{u\in H_{\Gamma_{0}}^{1}(\Omega) ;\; I(u) >0 \right\} \cup \{ 0\} \;
\mbox{ and } \;
\mathcal{N}^{-}=\left\{ u \in H_{\Gamma_{0}}^{1}(\Omega) ;I(u) < 0 \right\} .
\end{equation*}
The \textit{stable} set $\mathcal{W}$ and \textit{unstable} set $\mathcal{U}$ are defined respectively as:
\begin{equation*}
\mathcal{W=}\left\{ u\in H_{\Gamma_{0}}^{1}(\Omega) ;J(u) \leq d\right\} \cap \mathcal{N}^{+} \;
\mbox{ and } \;
\mathcal{U=}\left\{ u \in H_{\Gamma_{0}}^{1}(\Omega);J(u) \leq d \right\} \cap \mathcal{N}^{-}.
\end{equation*}

It is readily seen that the potential depth $d$ is also characterized by (see \cite{GS06})
\begin{equation}\label{d_Nihari}
d= \min_{u\in \mathcal{N}} J\left( u\right) .
\end{equation}
As it was remarked by Gazzola and Squassina in \cite{GS06}, this alternative characterization of $d$ shows that
\begin{equation} \label{altd}
\beta = \mbox{dist}(0,\mathcal{N}) =\min_{u\in \mathcal{N}} \Vert \nabla u\Vert_2 = \sqrt{\frac{ 2 d p }{p-2}} > 0.
\end{equation}
In Lemma \ref{lemme1}, we would like to prove that if the initial
datum $u_0$ is in the set $\mathcal{N}^{+}$ and if the initial energy $E(0)$ is not large
(we will precise exactly how large may be the initial energy), then $u(t)$ stays in $\mathcal{N}^{+}$, for each
$t\in [0,T)$, where $u(t)$ is the solution of (\ref{ondes}) obtained in Theorem \ref{existence}.

For this purpose, as in \cite{GS06,V99}, we denote by $C_{\ast}$ the best constant in the Poincar\'{e}-Sobolev
embedding $H_{\Gamma_{0}}^{1}(\Omega) \hookrightarrow L^{p}(\Omega)$ defined by:
\begin{equation}\label{sobolev}
C_{\ast}^{-1} = \inf\left\{\Vert \nabla u \Vert_2 : u \in  H_{\Gamma_{0}}^{1}(\Omega), \Vert u\Vert_p = 1 \right\}.
\end{equation}
Let us denote the Sobolev critical exponent:
$$
\bar{p}=\left\{
\begin{array}{cl}
\displaystyle \frac{2 N }{N-2} \;,&\mbox{ if } N \geq 3 \\
+\infty \;,&\mbox{ if } N=1, 2
\end{array}
\right. .
$$
Let us remark (as in \cite{GS06,V99}) that if  $p < \bar{p}$ the previous embedding is compact and
the infimum in (\ref{sobolev}) (as well as in (\ref{potentialwell})) is attained.
In such case (see, e.g. \cite[Section 3]{PS75}), any mountain pass solution of the stationary problem is a minimizer for
(\ref{sobolev}) and $C_{\ast}$ is
related to its energy:
\begin{equation} \label{mountainpass}
d = \frac{p-2}{2 p}\; C_{\ast}^{-2 p/(p-2)}.
\end{equation}
Let us remark also that in the Theorem  \ref{existence}, we have supposed that $p < \bar{q}$ where $\bar{q}$ is
defined by (\ref{qbar}).
As $\bar{q} < \bar{p}$, we may use the above characterization of the potential well depth $d$.

We can now proceed in the global existence result investigation. For this sake, let us state three lemmas.

%%%%%%%%%%%%%%%%%%%%%%%%%%%%%%%%%%%%%%%%%%%%%%%%%%%%%%%%%%%%%%%%%%%%%%%%%%%%
% LEMMA 1
%%%%%%%%%%%%%%%%%%%%%%%%%%%%%%%%%%%%%%%%%%%%%%%%%%%%%%%%%%%%%%%%%%%%%%%%%%%%
\begin{lemma}\label{lemme1} Assume $2\leq p\leq \bar{q}$ and $\max\left( 2,\displaystyle \frac{\bar{q}}{\bar{q}+1-p} \right) \leq m \leq \bar{q}$. \\
Let $u_{0}\in \mathcal{N}^{+} \,,\ u_{0} \neq 0$ and $u_{1}\in L^{2}(\Omega) $. Moreover, assume that $E(0) < d$.
Let us define $u$ the solution of problem (\ref{ondes}) in the sense of the Defintion \ref{generalised}. Then $u(t,.) \in \mathcal{N}^{+}$ for each 
$t\in [ 0,T_{max}).$
\end{lemma}
\begin{remark} \rm  Let us remark, that if there exists $\overline{t} \in [0,T_{max})$ such that
$$E(\overline{t}) < d \quad  and   \quad u(\overline{t}) \in  \mathcal{N}^{+}$$
the  same result stays true. It is the reason why we choose $\overline{t} = 0$.

Moreover, one can easily see that, from (\ref{mountainpass}), the condition $E(0) < d$ is equivalent to the inequality:
\begin{equation} \label{initial}
C_{\ast}^{p}\left(\frac{2p}{p-2} E(0) \right)^{\frac{p-2}{2}} < 1
\end{equation}
This last inequality will be used in the remaining proofs.
\end{remark}
%%%%%%%%%%%%%%%%%%%%%%%%%%%%%%%%%%%%%%%%%%%%%%%%%%%%%%%%%%%%%%%%%%%%%%%%%%%%
% PROOF OF LEMMA 1
%%%%%%%%%%%%%%%%%%%%%%%%%%%%%%%%%%%%%%%%%%%%%%%%%%%%%%%%%%%%%%%%%%%%%%%%%%%%
\begin{proof}
Since $I(u_{0}) >0$,  then by continuity, there exists $T_{\ast}\leq T_{max}$ such that
$I(u(t,.)) \geq 0,$ for all $t\in [0,T_{\ast })$. Since we have the relation:
\begin{equation}\label{Function_J}
J(t)   =  \frac{p-2}{2p} \Vert \nabla u\Vert_{2}^{2}+ \frac{1}{p} I(t), \quad \forall t\in [0,T_{\ast })\nonumber\\
\end{equation}
we easily obtain :
\begin{equation*}
J(t) \geq \frac{p-2}{2p} \Vert \nabla u\Vert_{2}^{2}, \quad \forall t\in [0,T_{\ast }).
\end{equation*}
Hence we have:
\begin{equation*}
\Vert \nabla u \Vert_{2}^{2} \leq  \frac{2p}{p-2}J(t), \quad \forall t\in [0,T_{\ast }).
\end{equation*}
From (\ref{Energy_J}) and (\ref{Energy_E}), we obviously have $J(t) \leq E(t),\,\forall t \in [0,T_{\ast})$. Thus we obtain:
\begin{equation*}
\Vert \nabla u \Vert_{2}^{2} \leq  \frac{2p}{p-2}E(t),\qquad \forall t \in [0,T_{\ast}).
\end{equation*}
Since $E$ is a decreasing function of $t$, we finally have:
\begin{equation} \label{ineqE0}
\Vert \nabla u \Vert_{2}^{2} \leq \frac{2p}{p-2}E(0) ,\qquad\forall t\in [ 0,T_{\ast }) \;.
\end{equation}
By definition of $C_{\ast}$, we have:
\begin{equation*}
\Vert u \Vert_{p}^{p} \leq  C_{\ast }^{p} \Vert \nabla u \Vert_{2}^{p} =
                              C_{\ast }^{p} \Vert \nabla u \Vert_{2}^{p-2}\Vert \nabla u \Vert_{2}^{2} \ .
\end{equation*}
Using the inequality (\ref{ineqE0}), we deduce:
\begin{equation}\label{L_p_E_0_inequality}
\Vert u \Vert_{p}^{p} \leq C_{\ast }^{p} \left(\frac{2p}{p-2}E(0)\right)^{\frac{p-2}{2}}
\Vert \nabla u\Vert_{2}^{2}, \qquad \forall t \in [0,T_{\ast}).
\end{equation}
Now exploiting the inequality on the initial condition (\ref{initial}) we obtain:
\begin{equation*}
\Vert u \Vert_{p}^{p} < \Vert \nabla u\Vert_{2}^{2}, \qquad \forall t \in [0,T_{\ast}).
\end{equation*}
Hence $\Vert \nabla u \Vert_{2}^{2} - \Vert u \Vert_{p}^{p}>0, \; \forall t~\in~[ 0,T_{\ast })$. This shows that
$u(t,.) \in \mathcal{N}^{+}, \; \forall~t\in~[ 0,T_{\ast})$.
Since the energy $E$ is decreasing along trajectories, we have the following inequality:
\begin{equation*}
\lim_{t\rightarrow T_{\ast }}C_{\ast }^{p}\left[ \frac{2p}{p-2}E(t)\right]^{\frac{p-2}{2}}\leq C_{\ast }^{p} \left[ \frac{2p}{p-2}E(0) \right]^{\frac{p-2}{2}}<1,
\end{equation*}
Thus by repeating this procedure, $T_{\ast }$ is extended to $T_{max}$.
\end{proof}
%%%%%%%%%%%%%%%%%%%%%%%%%%%%%%%%%%%%%%%%%%%%%%%%%%%%%%%%%%%%%%%%%%%%%%%%%%%%
% LEMMA 2
%%%%%%%%%%%%%%%%%%%%%%%%%%%%%%%%%%%%%%%%%%%%%%%%%%%%%%%%%%%%%%%%%%%%%%%%%%%%
\begin{lemma}\label{lemme2}Assume $2\leq p\leq \bar{q}$ and
$\max\left( 2,\displaystyle \frac{\bar{q}}{\bar{q}+1-p} \right) \leq m \leq \bar{q}$. \\
Let $u_{0}\in \mathcal{N}^{+} \,,\ u_{0} \neq 0$ and $u_{1}\in L^{2}(\Omega) $. Moreover, assume that $E(0) < d$.
Then the solution of the problem (\ref{ondes}) in the sense of the Definition \ref{generalised} is global in time.
\end{lemma}
%%%%%%%%%%%%%%%%%%%%%%%%%%%%%%%%%%%%%%%%%%%%%%%%%%%%%%%%%%%%%%%%%%%%%%%%%%%%
% PROOF OF LEMMA 2
%%%%%%%%%%%%%%%%%%%%%%%%%%%%%%%%%%%%%%%%%%%%%%%%%%%%%%%%%%%%%%%%%%%%%%%%%%%%
\begin{proof} Since the map $t \mapsto E(t)$ is a non increasing function of time
 $t$, and using the relation (\ref{Function_J}), we have:
$$
E(0)  \geq E(t) =\frac{1}{2}\Vert u_{t} \Vert_{2}^{2}+\frac{1}{2}\Vert u_{t}\Vert_{2,\Gamma_{1}}^{2}+
         \frac{(p-2) }{2p}\Vert \nabla u \Vert_{2}^{2}+\frac{1}{p}I(t)\,,\qquad \forall t\in [0,T_{max}).
$$
By Lemma \ref{lemme1}, we know that $u(t,.)\in \mathcal{N}^{+}$ for all $t\in(0,T]$. Hence,
$$ E(0) \geq \frac{1}{2}\Vert u_{t}\Vert_{2}^{2}+
\frac{(p-2) }{2p} \Vert \nabla u \Vert_{2}^{2},\qquad \forall t\in [0,T_{max}).
$$
Thus, $\forall t \in [0,T_{max})\,,\, \mbox{ the norm } \Vert \nabla u \Vert_2  + \Vert u_t \Vert_2$ is uniformly bounded by a constant depending
only on $E(0)$ and $p$. Then by Definition \ref{Tmax}, the solution is global, that is  $T_{max} = \infty$.
\end{proof}

The following Lemma is crucial in the proof of our result. A similar one (but for a different problem) was introduced in \cite{E03}.
\begin{lemma}\label{stable_unstable}
For every solution of (\ref{ondes}), given by Theorem \ref{existence}, only one of the following assumption holds:
\begin{enumerate}
\item[(i)] if there exists some $\overline{t}\geq 0 \mbox{ such that } u(\overline{t}) \in \mathcal{W} \mbox{ and } E(\overline{t})<d $, then  $u(t) \in \mathcal{W}  \mbox{ and }  E(t)<d,\,\forall t \geq \overline{t}$.
\item[(ii)]  if there exists some $\overline{t} \geq 0 \mbox{ such that } u(\overline{t}) \in \mathcal{U}  \mbox{ and } E(\overline{t})<d $, then  $u(t) \in \mathcal{U}  \mbox{ and }  E(t)<d,\, \forall t \geq \overline{t}$.
\item[(iii)]  $ E(t) \geq d,\,\forall t \geq 0$.
\end{enumerate}
\end{lemma}
%%%%%%%%%%%%%%%%%%%%%%%%ùù
%Proof of Lemma 2.1
%%%%%%%%%%%%%%%%%%%%%%%%%
\begin{proof}
Without loss of generality, we may assume that $\overline{t}=0$ and all along the paper, we suppose that $u_0\neq0$. 

Let us first prove (i).
Indeed, exploiting inequality (\ref{derivE}), we deduce that the energy functional is  a non-increasing function and
consequently, $E(t)<d$, for all $t\in[ 0,T_{max})$. Therefore (\ref{Energy_E}) implies that
 $J(t)<d$ for all $t\in[ 0,T_{max})$. This together with Lemma \ref{lemme1} gives (i).
 
Secondly, let us prove (ii). Let $u_0\in \mathcal{U}$ such that $E(0)<d$. Then (\ref{derivE}) 
implies that 
$$E(t)\leq E(0)<d, \qquad\forall t\in [0,T_{max}).$$ 
Next, let us assume by contradiction that there exists $\hat{t}\in [0,T_{max})$ such
 that $u(\hat{t})\notin\mathcal{U}$ and by continuity $I(u(\hat{t}))=0$. This implies that $u(\hat{t})\in \mathcal{N}.$
 Now using (\ref{d_Nihari}), we get $J(u(\hat{t}))\geq d$. This cannot be true since $J(u(t))<d$, for all $t\in[ 0,T_{max})$.
 Consequently, (ii) holds. 
 
The assertion (iii) is always true if (i) and (ii) are false. This completes the proof of Lemma \ref{stable_unstable}. 
\end{proof}
We can now state the asymptotic behavior of the solution of problem (\ref{ondes}).
%%%%%%%%%%%%%%%%%%%%%%%%%%%%%%%%%%%%%%%%%%%%%%%%%%%%%%%%%%%%%%%%%%%%%%%%%%%%
% THEOREM : EXPONENTIAL DECAY
%%%%%%%%%%%%%%%%%%%%%%%%%%%%%%%%%%%%%%%%%%%%%%%%%%%%%%%%%%%%%%%%%%%%%%%%%%%%
\begin{theorem} \label{exponential} Assume $2\leq p\leq \bar{q}$ and
$\max\left( 2,\displaystyle \frac{\bar{q}}{\bar{q}+1-p} \right) \leq m \leq \bar{q}$. Let
$u_{0}\in \mathcal{N}^{+}$ and $u_{1}\in L^{2}(\Omega) $. Moreover, assume that $E(0) < d$.
Then there exist two positive constants $\widehat{C}$ and $\xi $ independent of  $t$ such that:
\begin{eqnarray*}
0 < E(t) \leq \widehat{C}e^{-\xi t},\qquad \forall \, t\geq 0.
\end{eqnarray*}
\end{theorem}
%%%%%%%%%%%%%%%%%%%%%%%%%%%%%%%%%%%%%%%%%%%%%%%%%%%%%%%%%%%%%%%%%%%%%%%%%%%%
% Remark
%%%%%%%%%%%%%%%%%%%%%%%%%%%%%%%%%%%%%%%%%%%%%%%%%%%%%%%%%%%%%%%%%%%%%%%%%%%%
\begin{remark} \rm
Let us remark that these inequalities imply that there exist positive constants $K$ and $\zeta$ independent of  $t$ such that:
\begin{eqnarray*}
\Vert \nabla u(t)\Vert_2^2 + \Vert u_t(t) \Vert_2^2 \leq K e^{-\zeta t},\qquad \forall \, t\geq 0.
\end{eqnarray*}
Thus, this result improves the decay rate of Gazzola and Squassina \cite[Theorem 3.8]{GS06} (although the problem investigated by the two authors is slightly different), in which they
showed only the polynomial decay of the wave equation with strong damping and Dirichlet boundary conditions on the whole boundary of the domain. Here we show that for any initial data satisfying $u_{0}\in \mathcal{N}^{+}$
and $u_{1}\in L^{2}(\Omega) $ and verify the inequality (\ref{initial}), the solution can decay
faster than $1/t$, in fact with an exponential rate, even in the case $m>2$.

Also, by adapting the following proof in the spirit of the work done by  Gazzola and Squassina in \cite{GS06}, we can show an exponential decay rate 
even in the absence of the strong damping ($\alpha=0$) and $m=2$.
\end{remark}
%%%%%%%%%%%%%%%%%%%%%%%%%%%%%%%%%%%%%%%%%%%%%%%%%%%%%%%%%%%%%%%%%%%%%%%%%%%%
% PROOF OF THEOREM
%%%%%%%%%%%%%%%%%%%%%%%%%%%%%%%%%%%%%%%%%%%%%%%%%%%%%%%%%%%%%%%%%%%%%%%%%%%%
\begin{proof}
Since $u_{0}\in \mathcal{N}^{+}$ and $E(0) < d$, by Lemma \ref{lemme1} and Lemma \ref{lemme2}, we already have $u(t)\in \mathcal{N}^{+}$ for all $t \geq 0$. So we firstly get:
\begin{eqnarray*}
0 < E(t), \qquad \forall \, t\geq 0.
\end{eqnarray*}
The proof of the other inequality relies on the construction of a Lyapunov functional by performing a suitable
modification of  the energy. To this end, for $\varepsilon >0$, to be chosen later, we define for $ u \in  \mathcal{N}^{+}$,
\begin{equation}\label{energy_L}
\forall t \geq 0  \;,\;L(t) =E(t) +\varepsilon \int_{\Omega}u_{t} u dx + \varepsilon \int_{\Gamma_{1}} u u_{t} d\sigma +
\frac{\varepsilon \alpha }{2} \Vert \nabla u\Vert_{2}^{2}.
\end{equation}
%%%%%%%%%%%%% NEW %%%%%%%%%%%
Let us see that we have: for all $t \geq 0$
\begin{equation*}
\left\vert L(t)-E(t) \right\vert=\left| \varepsilon \int_{\Omega}u_{t} u dx + \varepsilon \int_{\Gamma_{1}} u u_{t} d\sigma + \frac{\varepsilon \alpha }{2} \Vert \nabla u\Vert_{2}^{2} \right|.
\end{equation*}
Since we have proved in Lemma \ref{lemme1} and Lemma  \ref{lemme2} that for all $t \geq 0$
$I(t)>0 \mbox{ and } \Vert \nabla u \Vert_2  + \Vert u_t \Vert_2 $ is uniformly bounded by a constant depending
only on $E(0)$ and $p$, using Young's inequalities on the two integral terms and then Poincar\'{e}'s inequality,  there exists a constant $C > 0$ such that:
\begin{equation*}
\left|\varepsilon \int_{\Omega}u_{t} u dx + \varepsilon \int_{\Gamma_{1}} u u_{t} d\sigma +
\frac{\varepsilon \alpha }{2} \Vert \nabla u\Vert_{2}^{2} \right|\leq C\varepsilon E(t).
\end{equation*}
Consequently, from the above two inequalities, we have
\begin{equation*}
(1-C\varepsilon)E(t)\leq L(t)\leq (1+C\varepsilon)E(t), \qquad \forall t \geq 0.
\end{equation*}
It is clear that for $\varepsilon$ sufficiently small, we can find two positive constants $\beta_{1}$ and $\beta_{2}$ such that
\begin{equation} \label{equivLE}
\beta_{1}E(t) \leq L(t) \leq \beta_{2}E(t), \qquad \forall t \geq 0.
\end{equation}
%%%%%%%%%%%%% FIN NEW %%%%%%%%%%%
By taking the time derivative of the function $L$ defined above in equation (\ref{energy_L}),
 using problem (\ref{ondes}) and formula (\ref{derivE}),
and performing several integration by parts, we get:
\begin{eqnarray}
\frac{dL(t) }{dt} &=&-\alpha \Vert \nabla u_{t}\Vert_{2}^{2}-r\Vert u_{t}\Vert_{m,\Gamma_{1}}^{m}+\varepsilon
\Vert u_{t}\Vert_{2}^{2}-\varepsilon \Vert \nabla
u\Vert_{2}^{2} \nonumber\\
&&+\varepsilon \Vert u\Vert_{p}^{p}+\varepsilon \Vert
u_{t}\Vert_{2,\Gamma_{1}}^{2}-\varepsilon r\int_{\Gamma_{1}}\vert u_{t}\vert ^{m-2}u_{t}ud\sigma \label{dLdt}.
\end{eqnarray}%
Now, we estimate the last term in the right hand side of (\ref{dLdt}) as follows.\\
By using Young's inequality, we obtain, for any $\delta>0$
\begin{equation}\label{int_gamma1}
\left\vert \int_{\Gamma_{1}}\vert u_{t}\vert ^{m-2}u_{t}ud\sigma \right\vert \leq \frac{\delta ^{-m}}{m}\Vert u\Vert_{m,\Gamma_{1}}^{m}+
\frac{m-1}{m} \delta ^{m/\left( m-1\right) }\Vert u_{t}\Vert_{m,\Gamma_{1}}^{m}.
\end{equation}
The trace inequality  implies that:
$$
\Vert u\Vert_{m,\Gamma_{1}}^{m} \leq C\Vert \nabla u\Vert_{2}^{m} \ ,
$$
where $C$ here and in the sequel denotes a generic positive constant which
might change from line to line.
Since the inequality (\ref{ineqE0}) holds, we have
\begin{equation}\label{norm_u_m_gamma1}
\Vert u\Vert_{m,\Gamma_{1}}^{m} \leq C \left(\frac{ 2 \, p \, E(0)}{p-2}\right)^{\frac{m-2}{2}} \Vert \nabla u\Vert_{2}^{2}.
\end{equation}
Inserting the two inequalities  (\ref{int_gamma1}) and (\ref{norm_u_m_gamma1}) in (\ref{dLdt}) and using (\ref{L_p_E_0_inequality}), we have:
\begin{eqnarray}
\frac{dL(t)}{dt} &\leq &-\alpha \Vert \nabla
u_{t}\Vert_{2}^{2}+r
\left( \varepsilon \frac{m-1}{m}\delta^{m/(m-2)}-1\right) \Vert u_{t}\Vert_{m,\Gamma_{1}}^{m}\nonumber
\\
&&+\varepsilon \Vert u_{t}\Vert_{2}^{2}+\varepsilon \Vert u_{t}\Vert_{2,\Gamma_{1}}^{2} \label{estimdLdt1}\\
&&+\varepsilon \left( \frac{r \delta^{-m}}{m} C \left(\frac{ 2 \, p \, E(0)}{p-2}\right)^{\frac{m-2}{2}}
+\underset{<0}{\underbrace{C_{\ast }^{p}\left( \frac{2p}{\left(p-2\right) }E(0) \right) ^{\frac{p-2}{2}}-1}}
\right) \Vert\nabla u\Vert_{2}^{2}.\nonumber
\end{eqnarray}
From (\ref{initial}), we have
$$
\displaystyle C_{\ast }^{p}\left( \frac{2p} {\left(p-2\right)}
E(0) \right)^{\frac{p-2}{2}}-1 < 0.
$$
Now, let us choose $\delta$ large enough such that:
$$
\left(
\frac{r \delta^{-m}}{m} C \left(\frac{ 2 \, p \, E(0)}{p-2}\right)^{\frac{m-2}{2}}
+ C_{\ast}^{p}\left(\frac{2p}{p-2}E(0)\right)^{\frac{p-2}{2}}-1
\right) < 0.
$$
Once $\delta$ is fixed, we fix $\varepsilon$ small enough such that:
$$
\left( \varepsilon \frac{m-1}{m}\delta^{m/(m-2)}-1\right) < 0.
$$
From (\ref{estimdLdt1}), we may find $\eta > 0 $, which depends only on $\delta$, such that:
\begin{equation*}
\frac{dL(t) }{dt}\leq -\alpha \Vert \nabla u_{t}\Vert_{2}^{2}+\varepsilon \Vert u_{t}\Vert_{2}^{2}+
\varepsilon \Vert u_{t}\Vert_{2,\Gamma_{1}}^{2}-\varepsilon \eta \Vert \nabla u\Vert_{2}^{2}.
\end{equation*}
Consequently, using the definition of the energy (\ref{Energy_E}), for any positive constant $M$, which will be chosen below,  we obtain:
\begin{eqnarray}
\frac{dL(t)}{dt} &\leq &-M\varepsilon E(t) +\varepsilon \left( 1+\frac{M}{2}\right)\Vert u_{t}\Vert_{2}^{2}-
\alpha \Vert \nabla u_{t}\Vert_{2}^{2}  \nonumber \\
&&+\left( \frac{M\varepsilon }{2}+\varepsilon \right) \Vert u_{t}\Vert_{2,\Gamma_{1}}^{2}+
\varepsilon \left( \frac{M}{2}-\eta \right) \Vert \nabla u\Vert_{2}^{2}. \label{estimdLdt2}
\end{eqnarray}
By using the Poincar\'{e} inequality and the trace inequality
\begin{eqnarray*}
\Vert u_{t}\Vert_{2}^{2} &\leq &C\Vert \nabla u_{t}\Vert_{2}^{2} \\
\Vert u_{t}\Vert_{2,\Gamma_{1}}^{2} &\leq &C\Vert \nabla u_{t}\Vert_{2}^{2},
\end{eqnarray*}
choosing again $\varepsilon$ small enough and $M\leq 2\eta$, from (\ref{estimdLdt2}), we have:
\begin{equation*}
\frac{dL(t) }{dt}\leq -M\varepsilon E(t) ,\qquad \forall t \geq 0.
\end{equation*}%
On the other hand, by virtue of (\ref{equivLE}), setting $\xi =M\varepsilon /\beta_{2}$, the last inequality becomes:
\begin{equation}\label{diffineq}
\frac{dL(t)}{dt} \leq -\xi L(t) \;,\qquad \forall t\geq0.
\end{equation}
Integrating the previous differential inequality (\ref{diffineq}) between $0$ and $t$ gives the following estimate for the
function $L$:
\begin{equation*}
L(t) \leq Ce^{-\xi t} \;,\qquad \forall t\geq0.
\end{equation*}%
Consequently, by using  (\ref{equivLE}) once again, we conclude
\begin{equation*}
E(t) \leq \widehat{C}e^{-\xi t} \;,\qquad \forall t\geq0 .
\end{equation*}
This completes the proof of Theorem \ref{exponential}.
\end{proof}
\begin{remark} \rm
In \cite{GS08}, we have proved the following result:
\begin{theorem}
Assume $2\leq p\leq \bar{q}$ and
$m < p $. Let $u_{0}\in H_{\Gamma_{0}}^1(\Omega)$ and $u_{1}\in L^{2}(\Omega) $.

Suppose that
$$
E(0) <d \mbox{ and } \left\Vert \nabla u_{0}\right\Vert _{2}>C_\ast^{-p/(p-2)}.
$$
Then the solution of problem (\ref{ondes}) growths exponentially in the $L^{p}$ norm.
\end{theorem}
The present result on the asymptotic stability completes the above result on the exponential growth since
when $u_{0}\in \mathcal{N}^{+}$, we have:  $\left\Vert \nabla u_{0}\right\Vert _{2} \leq C_\ast^{-p/(p-2)}$.

Indeed, since $d$ is the mountain pass level of the function $J$, we have $J(u_{0}) \leq d$.
From (\ref{Function_J}), this writes:
$$J(u_0)=\frac{p-2}{2p} \Vert \nabla u_0\Vert_{2}^{2}+ \frac{1}{p} I(0) \leq d $$
Since $u_0 \in \mathcal{N}^{+}$, we have $I(0)>0$ and consequently,
$$\frac{p-2}{2p} \Vert \nabla u_0\Vert_{2}^{2} \leq d.$$
Using identity (\ref{mountainpass}), we get finally
$\left\Vert \nabla u_{0}\right\Vert _{2} \leq C_\ast^{-p/p-2}$.
\end{remark}
%%%%%%%%%%%%%%%%%%%%%%%%%%%%%%%%%%%%%%%%%%%%%%%%%%%%%%%%%%%%%%%%%%%%%%%%%%%%%%%%%%%%%%%%%
\section{Blow up}
%%%%%%%%%%%%%%%%%%%%%%%%%%%%%%%%%%%%%%%%%%%%%%%%%%%%%%%%%%%%%%%%%%%%%%%%%%%%%%%%%%%%%%%%%
In this section we consider the problem (\ref{ondes}) in the linear boundary damping case (i.e. $m=2$) and
we show that if for some $\overline{t} \in [0,T_{max})\,,\, u(\overline{t}) \in \mathcal{U}$ and $E(\overline{t}) \leq d$ then the solution
of (\ref{ondes}) blows up in finite time. Our result reads as follows:
\begin{theorem}\label{blowup}
Assume $2\leq p\leq \bar{q}$ and $m=2$. Let $u$ be the solution of (\ref{ondes}) on $[0,T_{max})$.
Then $T_{max }<\infty$ if and only if there exists $\overline{t} \in [0,T_{max})$ such that:
\begin{equation}\label{blowcondition}
u(\overline{t})\in \mathcal{U}\ \mbox{  and } \ E(\overline{t}) \leq d.
\end{equation}
\end{theorem}

%%%%%%%%%%%%%%%%%%%%%%%%%%%%%%%%%%%%%%%%%%%%%%%%%%%%%%%%%%%%%%%%%%%%%%%%%%%%
% PROOF OF BLOW UP
%%%%%%%%%%%%%%%%%%%%%%%%%%%%%%%%%%%%%%%%%%%%%%%%%%%%%%%%%%%%%%%%%%%%%%%%%%%%
\begin{proof}Without loss of generality, we may assume that $\overline{t}=0$.

Let us suppose that $u(0)\in \mathcal{U}\mbox{ and }  E(0) \leq d$.
We will prove that $T_{max }<\infty$ by contradiction. We will suppose that the solution is global ``in time''
and we will use the concavity argument due to Levine
\cite{L74_1,L74_2} where the basic idea of this method is to construct a positive functional $\theta(t)$ of the solution
and show that for some $\gamma >0$, the function $\theta^{-\gamma}(t) $ is a positive concave function of $t$.
Thus it will exist $T^*$ such that $\displaystyle\lim_{t\rightarrow T^*} \theta^{-\gamma}(t) = 0$.
From the construction of the function $\theta$, this will imply that:
$$
\lim_{\underset {t < T^*} {t \rightarrow T^*}} \Vert \nabla u \Vert_2  + \Vert u_t \Vert_2 = + \infty .
$$
In order to find such $\gamma$, we will verify that:
\begin{equation}\label{concav}
\frac{d^{2}\theta^{-\gamma }(t) }{dt^{2}}=-\gamma \theta^{-\gamma -2}(t)
\left[ \theta \theta^{^{\prime \prime}}-(1+\gamma) \theta ^{^{\prime 2}}(t) \right]
\leq 0 \; , \qquad \forall t \geq 0.
\end{equation}%
Thus it suffices to prove that $\theta(t)$ satisfies the differential inequality
\begin{equation}\label{difftheta}
\theta \theta^{^{\prime \prime }}-\left( 1+\gamma \right) \theta ^{^{\prime 2}}(t) \geq 0\; , \qquad \forall t \geq 0.
\end{equation}
From Lemma \ref{stable_unstable}, we firstly have:
$$ E(t) \leq  d \mbox{ and }
u(t) \in \mathcal{U}, \qquad \forall t\in [ 0,T_{max}).$$
It is clear that $\mathcal{N}$ can be seen as a set which separate the two 
sets $\mathcal{N}^+$ and $\mathcal{N}^-$ in $H_{\Gamma_0}^{1}$. 

From the definition (\ref{potentialwell})
of the potential well  depth $d$, for $u\in H_{\Gamma_{0}}^{1}(\Omega)\backslash \{0\}$, we have:
\begin{equation}
d\leq \sup_{\lambda \geq 0}J(\lambda u) =\frac{p-2}{2p}\left(\frac{\Vert \nabla u\Vert_{2}^{2p}}{\Vert u\Vert_{p}^{2p}}\right)^{\frac{1}{(p-2)}}.  \label{d}
\end{equation}
On the other hand, since $\forall t\in [0,T_{max})\,,\, u \in \mathcal{N}^{-}$, we have:
$$\forall t\in [0,T_{max})\,,\,I(t) <0 \quad . $$
This inequality gives naturally $\forall t\in [0,T_{max})\,,\,\Vert\nabla u\Vert_{2}^{2}<\Vert u\Vert_{p}^{p}$. 
Therefore, using this last inequality, the inequality (\ref{d}) becomes:
\begin{equation*}
\forall t\in [0,T_{max})\,,\, d<\frac{p-2}{2p}\left\Vert \nabla u\right\Vert_{2}^{2},
\end{equation*}
which will be used as:
\begin{equation} \label{estim_grad_u}
\frac{2dp}{p-2}<\Vert \nabla u(t) \Vert_{2}^{2} \, , \qquad \forall t\in \left[ 0,T_{max }\right) \, .
\end{equation}
Assume by contradiction that the solution $u$ is global ``in time''. Then for any $T > 0$, let us define the functional
$\theta $ as follows
\begin{eqnarray}
\theta (t) &=&\Vert u(t) \Vert_{2}^{2} + \Vert u(t) \Vert_{2,\Gamma_{1}}^{2} +
\alpha \int_{0}^{t}\Vert \nabla u\left( s\right) \Vert_{2}^{2}ds
+ r\int_{0}^{t}\Vert u\left( s\right) \Vert_{2,\Gamma_{1}}^{2}ds  \nonumber \\
&& + (T-t) \left[ \alpha \Vert \nabla u_{0}\Vert_{2}^{2}+r\Vert u_{0}\Vert_{2,\Gamma_{1}}^{2}\right] ,
\qquad \forall t\in [ 0,T). \label{deftheta}
\end{eqnarray}%
Taking the time derivative of (\ref{deftheta}) we have:
\begin{eqnarray}
\theta^{\prime}(t) &=& 2\int_{\Omega }u_{t}u dx + 2\int_{\Gamma_{1}}u_{t}u d\sigma +
2\alpha \int_{0}^{t}\int_{\Omega }\nabla u\nabla u_{t}dxds  \nonumber \\
&&+2r\int_{0}^{t}\int_{\Gamma_{1}}u_{t}ud\sigma ds. \label{thetaprime}
\end{eqnarray}
Replacing $u_{tt}$ by its expression given by problem (\ref{ondes}) and using Green's formula (see  \cite{LM68}), the function
$\theta'$ is differentiable and we have:
\begin{equation*}
\theta^{\prime \prime}(t) =2 \left[ \Vert u_{t}(t) \Vert_{2}^{2} - \Vert \nabla u(t) \Vert_{2}^{2} + \Vert u\Vert_{p}^{p} + \Vert u_{t}(t) \Vert_{2,\Gamma_{1}}^{2}\right] .
\end{equation*}
Therefore, using the definition of $\theta$ given by (\ref{deftheta}), we can easily see that:
\begin{eqnarray}
\theta (t) \theta^{^{\prime \prime }}(t) &-&\frac{p+2}{4}\theta ^{^{\prime }}(t) ^{2}
=2\theta (t) \left[ \Vert u_{t}(t) \Vert_{2}^{2}-\Vert \nabla u(t) \Vert _{2}^{2}+\Vert u\Vert_{p}^{p}+\Vert u_{t}(t)
\Vert_{2,\Gamma_{1}}^{2}\right]  \nonumber \\
&&-\left( p+2\right) \biggl[ \theta (t) -
\left( T-t\right) \Bigl[\alpha \Vert \nabla u_{0}\Vert_{2}^{2}+ r\Vert u_{0}\Vert_{2,\Gamma_{1}}^{2}\Bigr] \biggr]
\label{eqtheta}\\
&&\times \left[ \Vert u_{t}(t) \Vert
_{2}^{2}+\Vert u_{t}(t) \Vert_{2,\Gamma
_{1}}^{2}+\alpha \int_{0}^{t}\Vert \nabla u_{t}(t)
\Vert_{2}^{2}ds+r\int_{0}^{t}\Vert u_{t}(t)
\Vert_{2,\Gamma_{1}}^{2}ds\right]  \nonumber \\
&&+\left( p+2\right) \eta (t)  \nonumber
\end{eqnarray}
where the function $\eta$ is defined by:
\begin{eqnarray}
&&\eta(t) = \left[ \Vert u(t) \Vert_{2}^{2}+\Vert u(t) \Vert_{2,\Gamma_{1}}^{2}+\alpha \int_{0}^{t}\Vert \nabla u(t)
\Vert_{2}^{2}ds+r\int_{0}^{t}\Vert u(t) \Vert_{2,\Gamma_{1}}^{2}ds\right]  \nonumber \\
&&\times \left[ \Vert u_{t}(t) \Vert_{2}^{2}+\Vert u_{t}(t) \Vert_{2,\Gamma_{1}}^{2}+
\alpha \int_{0}^{t}\Vert \nabla u_{t}(t)\Vert_{2}^{2}ds+r\int_{0}^{t}\Vert u_{t}(t) \Vert_{2,\Gamma_{1}}^{2}ds\right]
\label{eta}  \\
&&-\left[ \int_{\Omega }u_{t}udx+\int_{\Gamma_{1}}u_{t}ud\sigma +
\alpha \int_{0}^{t}\int_{\Omega }\nabla u\nabla u_{t}dxds+r\int_{0}^{t}\int_{\Gamma_{1}}u_{t}ud\sigma ds\right]^{2}.\nonumber
\end{eqnarray}
Our purpose now is to show that the right hand side of the equality (\ref{eqtheta}) is non negative.
Let us firstly show that $\eta(t)\geq 0$ for every $t\in [0,T)$.
To do this, we estimate all the terms in the third line of (\ref{eta}) making use of Cauchy-Schwarz inequality, and
compare the results with the terms in the first and second line in (\ref{eta}).
For instance, when we develop the square term in the inequality (\ref{eta}), we estimate the terms as follows:
\begin{eqnarray*}
\left( \int_{\Omega }u_{t}udx\right) ^{2} &\leq &\Vert u(t) \Vert_{2}^{2}\Vert u_{t}(t) \Vert_{2}^{2}\; \mbox{ and }\\
2\int_{\Omega }u_{t} u dx\int_{\Gamma_{1}}u_{t} u d\sigma &\leq &\Vert u(t) \Vert_{2,\Gamma_{1}}^{2}
\Vert u_{t}(t) \Vert_{2}^{2}+\Vert u_{t}(t) \Vert_{2,\Gamma_{1}}^{2}\Vert u(t) \Vert_{2}^{2}.
\end{eqnarray*}
Also, the following estimate holds:
\begin{eqnarray*}
2\alpha \int_{0}^{t}\int_{\Omega }\nabla u\nabla u_{t}dxds\int_{\Omega}u_{t}u dx &\leq&
\alpha \Vert u_{t}(t) \Vert_{2}^{2}\int_{0}^{t}\Vert \nabla u\left( s\right) \Vert_{2}^{2}ds \\
&~& +\alpha \Vert u(t) \Vert_{2}^{2}\int_{0}^{t}\Vert \nabla u_{t}\left( s\right) \Vert_{2}^{2}ds.
\end{eqnarray*}%
By carrying ``carefully'' all computations based on the same estimates as above, we finally obtain
$$\eta (t) \geq 0, \qquad \forall t\in [0,T).$$
Consequently, the equality (\ref{eqtheta}) becomes
\begin{equation*}
\theta (t) \theta ^{^{\prime \prime }}(t) -\frac{p+2}{4}\theta ^{^{\prime }}(t) ^{2} \geq
\theta (t) \zeta (t),\qquad \forall t\in [0,T).
\end{equation*}
where
\begin{eqnarray*}
\zeta (t) &=&2\left[ \Vert u_{t}(t) \Vert_{2}^{2} - \Vert \nabla u(t) \Vert_{2}^{2} +
\Vert u\Vert_{p}^{p} + \Vert u_{t}(t) \Vert_{2,\Gamma_{1}}^{2}\right] \\
&&-(p+2) \Bigl\{ \Vert u_{t}(t) \Vert_{2}^{2}+\Vert u_{t}(t) \Vert_{2,\Gamma_{1}}^{2}\\
&&+\alpha \int_{0}^{t}\Vert \nabla u_{t}(t) \Vert_{2}^{2}ds +
r\int_{0}^{t}\Vert u_{t}(t)\Vert_{2,\Gamma_{1}}^{2}ds\Bigr\}  \ .\label{zeta}\\
\end{eqnarray*}%
Let us remark that
\begin{eqnarray*}
\zeta (t) &= & -2pE(t) +(p-2) \Vert \nabla u(t)\Vert_{2}^{2} - (p+2) \alpha
\int_{0}^{t}\Vert \nabla u_{t}\left( s\right) \Vert_{2}^{2}ds \\
&&-( p+2) r\int_{0}^{t}\Vert u_{t}(s)\Vert_{2,\Gamma_{1}}^{2}ds.
\end{eqnarray*}
From the equality (\ref{derivE}), we have:
\begin{equation}\label{EnergyE0}
E(t) +\alpha \int_{0}^{t}\Vert \nabla u_{t}(s) \Vert_{2}^{2}ds +
r\int_{0}^{t}\Vert u_{t}(s)\Vert_{2,\Gamma_{1}}^{2}ds=E(0), \qquad \forall t\in [0,T).
\end{equation}
Thus we can write:
\begin{eqnarray*}
\zeta (t) &=&-2pE(0) + (p-2) \Vert \nabla u(t) \Vert_{2}^{2} \\
&&+ (p-2) \alpha \int_{0}^{t}\Vert \nabla u_{t}(s) \Vert_{2}^{2}ds+
\left( p-2\right) r\int_{0}^{t}\Vert u_{t}(s) \Vert_{2,\Gamma_{1}}^{2}ds.
\end{eqnarray*}
Therefore, by using (\ref{estim_grad_u}) and since $E(0)\leq d$ we have:
\begin{eqnarray*}
\zeta (t) &>&2 p ( d-E(0)) + (p-2) \alpha \int_{0}^{t}\Vert \nabla u_{t}(s) \Vert_{2}^{2} ds+
(p-2) r \int_{0}^{t}\Vert u_{t}(s) \Vert_{2,\Gamma_{1}}^{2}ds \\
&\geq & (p-2) \alpha \int_{0}^{t}\Vert \nabla u_{t}(s) \Vert_{2}^{2}ds +
(p-2) r\int_{0}^{t}\Vert u_{t}(s) \Vert_{2,\Gamma_{1}}^{2}ds.
\end{eqnarray*}%
Hence, there exist  $t_{0}>0$ and $\delta >0$ such that
\begin{equation*}
\zeta (t) \geq \delta,\qquad \forall t \in [t_{0},T) \ .
\end{equation*}
Also, since $\theta (t) $ is continuous and positive, there exists $\rho >0$ such that
\begin{equation*}
\theta (t) \geq \rho,\qquad \forall t \in [t_{0},T)  \ .
\end{equation*}
Consequently,
\begin{equation*}
\theta (t) \theta^{\prime \prime}(t) -\frac{p+2}{4}\theta^{\prime}(t)^{2}\geq \rho \delta,\qquad \forall t \in  [t_{0},T)  \ .
\end{equation*}
Setting
$$ \gamma = \frac{ p - 2}{4} > 0,$$
the differential inequality (\ref{difftheta}) is verified on $[t_0,T)$. This proves that $\theta(t)^{-\gamma}$ reaches 0
in finite time, say as $t \rightarrow T^*$. Since $T^*$ is independent of the initial choice of $T$, we may assume that
 $T^* < T$. This tells us that:
$$ \lim_{t \rightarrow T^*} \theta(t) = + \infty.$$
From Poincar\'e's inequality and the continuity of the trace operator on $\Gamma_{1}$, by the equation  (\ref{deftheta}) defining $\theta$, this implies that:
$$
\lim_{\underset {t < T^*} {t \rightarrow T^*}} \Vert \nabla u \Vert_2 = + \infty.
$$
Thus we cannot suppose that the solution of (\ref{ondes}) with $m=2$ is global ``in time'', that is $T_{max} < \infty$.

Conversely, let us suppose that $T_{max}<\infty$. 
We want to show that there exists $\overline{t} \in [0,T_{max})$ such that:
\begin{equation*}
u(\overline{t})\in \mathcal{U} \mbox{ and }  E(\overline{t}) \leq d.
\end{equation*}
Notice first that, for every $0 < t < T_{max}$, by H\"{o}lder's inequality, there holds
\begin{equation*}
\displaystyle \bigint_{0}^{t}\Vert u_{t}(\tau)\Vert_{\ast}^{2} d\tau \geq \frac{1}{t}\left(\bigint_{0}^{t}\Vert u_{t}(\tau)\Vert_{\ast}d\tau\right)^{2} \ .
\end{equation*}
Thus since 
\begin{equation*}
\bigint_{0}^{t}\Vert u_{t}(\tau)\Vert_{\ast}d\tau\geq \left\Vert  \bigint_{0}^{t}u_{t}(\tau)d\tau \right\Vert_{\ast} 
\geq \Bigl\vert \Vert u(t)\Vert_{\ast} - \Vert u(0)\Vert_{\ast} \Bigl\vert\ ,
\end{equation*}
we have:
\begin{equation}\label{GazzSquass}
\bigint_{0}^{t}\Vert u_{t}(\tau)\Vert_{\ast}^{2} d\tau \geq \frac{1}{t}\Bigl(\Vert u(t)\Vert_{\ast} - \Vert u(0)\Vert_{\ast}\Bigl)^{2} \ .
\end{equation}
By the help of (\ref{EnergyE0}) and (\ref{GazzSquass}), we thus have:
\begin{equation}\label{Energy1}
E(t)\leq E(0)-\frac{1}{t}\Bigl(\Vert u(t)\Vert_{\ast}- \Vert u(0)\Vert_{\ast}\Bigl)^{2} . 
\end{equation}
To prove that the conditions (\ref{blowcondition}) are necessary, we will adapt the study of the dynamics of the waves equation performed by 
J. Esquivel-Avila in \cite{E03}.

We proceed by contradiction and we assume that for all $t\geq 0 \,,\, u(t) \notin \mathcal{U}$.
Then, by Lemma \ref{stable_unstable}, we have either:
\begin{enumerate}
\item[i)] $u( t) \in \mathcal{W}$ and $E(t) <d$ or
\item[ii)]  $E(t) \geq d$
\end{enumerate}
In the first case,  Lemma \ref{lemme2} implies that the solution is global in time. This is not possible since we assumed that $T_{max }<\infty $.
In the second case by using (\ref{Energy1}) we get
\begin{equation*}
\Bigl(\Vert u(t)\Vert_{\ast} - \Vert u(0)\Vert_{\ast}\Bigl)^2 \leq t \Bigl(E(0) - d\Bigl) \qquad \forall t \in [0,T_{max}).
\end{equation*}
Thus for any time $T  \in [0,T_{max})$ there exists a constant $C(T)$ such that
 \begin{equation}\label{Norm_star}
\Vert u(T)\Vert_{\ast }\leq C(T) .
\end{equation}
On the other hand, from Definition \ref{Tmax} and if $T_{max }<\infty $, then
\begin{equation*}
\lim_{t\rightarrow T_{max }}\left( \left\Vert u_{t}\left( t\right)
\right\Vert _{2}^{2}+\left\Vert \nabla u\left( t\right) \right\Vert_{2}^{2}\right) = + \infty \ .
\end{equation*}
Since the energy $E$ is decreasing along trajectories, we have the inequality:
\begin{equation*}
\frac{1}{2}\left( \Vert u_{t}(t)\Vert_{2,\Gamma_{1}}^{2}+\Vert u_{t}(t)\Vert_{2}^{2}+\Vert \nabla u(t) \Vert_{2}^{2}\right) \leq E(u(0)) +\frac{1}{p}\Vert u\Vert_{p}^{p},
\end{equation*}
This implies that
\begin{equation} \label{L_p_norm}
 \lim_{\underset {t < T_{max}} {t \rightarrow T_{max}}} \Vert u(t) \Vert_{p}=+\infty . 
\end{equation}
Combining (\ref{L_p_norm}) with the Poincar\'{e} inequality, we deduce that for every $M > E(0) $, there
exists $\hat{t}>0$ such that
\begin{equation}
M<\frac{p-2}{2p}\Vert u(\hat{t}) \Vert _{p}^{2}\leq \frac{p-2}{2p}\Vert \nabla u(\hat{t}) \Vert_{2}^{2}.  \label{M_inequality}
\end{equation}
This inequality is exactly the same as inequality (\ref{estim_grad_u}), where $d$ is replaced by $M$. 
Consequently following the proof of the sufficiency part of the theorem, by replacing $d$ by $M$ and 
by defining the same function $\theta $ but for $t \in [\hat{t},T_{max })$, we deduce that the solution
blows up in finite time $T^{\ast} \in (\hat{t},T_{max})$. Thus we obtain:
$$
\lim_{\underset {t < T_{\ast}} {t \rightarrow T_{\ast}}} \Vert \nabla u \Vert_2  = + \infty.
$$
and this contradicts (\ref{Norm_star}).
\end{proof}
\begin{remark} \rm
The term $f(u)=\vert u\vert^{p-2}u$ is clearly responsible for the blow up situation. It is often called the
``blow up term''.
Consequently when $f(u)=0$, or $f(u)=-\vert u\vert^{p-2}u$ any solution with arbitrary initial data is
global in time and the result of Theorem \ref{blowup} holds without condition (\ref{initial}).
\end{remark}
\begin{remark} \rm
 It's early well known (\cite{L74_1,L74_2}) that this blow up result appears for solutions with large initial data i.e.
$E(0)<0$. We note here that if $E(0)<0$, then the blow up conditions (\ref{blowcondition}) hold.
\end{remark}
\begin{ack}
The second author was supported by MIRA 2007 project of the R\'egion Rh\^one-Alpes.
This author wishes to thank Univ. de Savoie of Chamb\'ery for its kind hospitality.

Moreover, the two authors wish to thank deeply the referees for their useful remarks and their careful reading of the proofs presented in this paper.
\end{ack}
%%%%%%%%%%%%%%%%%%%%%%%%%%%
%
%%%%%%%%%%%%%%%%%%%%%%%%%%%%%
\bibliographystyle{siam}

\end{document}